\documentclass{article}

\usepackage{amssymb,amsmath,theorem,euscript}

\newcounter{sec}

\def\sm{\smallskip}

\newcounter{punct}[sec]

\def\punct{\refstepcounter{punct}{\arabic{sec}.\arabic{punct}.  }}

\def\COUNTERS{\addtocounter{sec}{1}
              \setcounter{punct}{0}
          \setcounter{equation}{0}
          \setcounter{theorem}{0}
                  }
\newtheorem{theorem}{Theorem}[sec]

\newtheorem{lemma}[theorem]{Lemma}

\begin{document}

 \def\ov{\overline}
\def\wt{\widetilde}
 \newcommand{\rk}{\mathop {\mathrm {rk}}\nolimits}
\newcommand{\Aut}{\mathop {\mathrm {Aut}}\nolimits}
\newcommand{\Out}{\mathop {\mathrm {Out}}\nolimits}
\renewcommand{\Re}{\mathop {\mathrm {Re}}\nolimits}
\def\Br{\mathrm {Br}}

\def\SL{\mathrm {SL}}
\def\SU{\mathrm {SU}}
\def\GL{\mathrm {GL}}
\def\U{\mathrm U}
\def\OO{\mathrm O}
 \def\Sp{\mathrm {Sp}}
 \def\SO{\mathrm {SO}}
\def\SOS{\mathrm {SO}^*}
 \def\Diff{\mathrm{Diff}}
 \def\Vect{\mathfrak{Vect}}
\def\PGL{\mathrm {PGL}}
\def\PU{\mathrm {PU}}
\def\PSL{\mathrm {PSL}}
\def\Symp{\mathrm{Symp}}
\def\End{\mathrm{End}}
\def\Mor{\mathrm{Mor}}
\def\Aut{\mathrm{Aut}}
 \def\PB{\mathrm{PB}}
 \def\cA{\mathcal A}
\def\cB{\mathcal B}
\def\cC{\mathcal C}
\def\cD{\mathcal D}
\def\cE{\mathcal E}
\def\cF{\mathcal F}
\def\cG{\mathcal G}
\def\cH{\mathcal H}
\def\cJ{\mathcal J}
\def\cI{\mathcal I}
\def\cK{\mathcal K}
 \def\cL{\mathcal L}
\def\cM{\mathcal M}
\def\cN{\mathcal N}
 \def\cO{\mathcal O}
\def\cP{\mathcal P}
\def\cQ{\mathcal Q}
\def\cR{\mathcal R}
\def\cS{\mathcal S}
\def\cT{\mathcal T}
\def\cU{\mathcal U}
\def\cV{\mathcal V}
 \def\cW{\mathcal W}
\def\cX{\mathcal X}
 \def\cY{\mathcal Y}
 \def\cZ{\mathcal Z}
\def\0{{\ov 0}}
 \def\1{{\ov 1}}
 \def\frA{\mathfrak A}
 \def\frB{\mathfrak B}
\def\frC{\mathfrak C}
\def\frD{\mathfrak D}
\def\frE{\mathfrak E}
\def\frF{\mathfrak F}
\def\frG{\mathfrak G}
\def\frH{\mathfrak H}
\def\frI{\mathfrak I}
 \def\frJ{\mathfrak J}
 \def\frK{\mathfrak K}
 \def\frL{\mathfrak L}
\def\frM{\mathfrak M}
 \def\frN{\mathfrak N} \def\frO{\mathfrak O} \def\frP{\mathfrak P} \def\frQ{\mathfrak Q} \def\frR{\mathfrak R}
 \def\frS{\mathfrak S} \def\frT{\mathfrak T} \def\frU{\mathfrak U} \def\frV{\mathfrak V} \def\frW{\mathfrak W}
 \def\frX{\mathfrak X} \def\frY{\mathfrak Y} \def\frZ{\mathfrak Z} \def\fra{\mathfrak a} \def\frb{\mathfrak b}
 \def\frc{\mathfrak c} \def\frd{\mathfrak d} \def\fre{\mathfrak e} \def\frf{\mathfrak f} \def\frg{\mathfrak g}
 \def\frh{\mathfrak h} \def\fri{\mathfrak i} \def\frj{\mathfrak j} \def\frk{\mathfrak k} \def\frl{\mathfrak l}
 \def\frm{\mathfrak m} \def\frn{\mathfrak n} \def\fro{\mathfrak o} \def\frp{\mathfrak p} \def\frq{\mathfrak q}
 \def\frr{\mathfrak r} \def\frs{\mathfrak s} \def\frt{\mathfrak t} \def\fru{\mathfrak u} \def\frv{\mathfrak v}
 \def\frw{\mathfrak w} \def\frx{\mathfrak x} \def\fry{\mathfrak y} \def\frz{\mathfrak z} \def\frsp{\mathfrak{sp}}
 \def\bfa{\mathbf a} \def\bfb{\mathbf b} \def\bfc{\mathbf c} \def\bfd{\mathbf d} \def\bfe{\mathbf e} \def\bff{\mathbf f}
 \def\bfg{\mathbf g} \def\bfh{\mathbf h} \def\bfi{\mathbf i} \def\bfj{\mathbf j} \def\bfk{\mathbf k} \def\bfl{\mathbf l}
 \def\bfm{\mathbf m} \def\bfn{\mathbf n} \def\bfo{\mathbf o} \def\bfp{\mathbf p} \def\bfq{\mathbf q} \def\bfr{\mathbf r}
 \def\bfs{\mathbf s} \def\bft{\mathbf t} \def\bfu{\mathbf u} \def\bfv{\mathbf v} \def\bfw{\mathbf w} \def\bfx{\mathbf x}
 \def\bfy{\mathbf y} \def\bfz{\mathbf z} \def\bfA{\mathbf A} \def\bfB{\mathbf B} \def\bfC{\mathbf C} \def\bfD{\mathbf D}
 \def\bfE{\mathbf E} \def\bfF{\mathbf F} \def\bfG{\mathbf G} \def\bfH{\mathbf H} \def\bfI{\mathbf I} \def\bfJ{\mathbf J}
 \def\bfK{\mathbf K} \def\bfL{\mathbf L} \def\bfM{\mathbf M} \def\bfN{\mathbf N} \def\bfO{\mathbf O} \def\bfP{\mathbf P}
 \def\bfQ{\mathbf Q} \def\bfR{\mathbf R} \def\bfS{\mathbf S} \def\bfT{\mathbf T} \def\bfU{\mathbf U} \def\bfV{\mathbf V}
 \def\bfW{\mathbf W} \def\bfX{\mathbf X} \def\bfY{\mathbf Y} \def\bfZ{\mathbf Z} \def\bfw{\mathbf w}
 \def\R {{\mathbb R }} \def\C {{\mathbb C }} \def\Z{{\mathbb Z}} \def\H{{\mathbb H}} \def\K{{\mathbb K}}
 \def\N{{\mathbb N}} \def\Q{{\mathbb Q}} \def\A{{\mathbb A}} \def\T{\mathbb T} \def\P{\mathbb P} \def\G{\mathbb G}
 \def\bbA{\mathbb A} \def\bbB{\mathbb B} \def\bbD{\mathbb D} \def\bbE{\mathbb E} \def\bbF{\mathbb F} \def\bbG{\mathbb G}
 \def\bbI{\mathbb I} \def\bbJ{\mathbb J} \def\bbL{\mathbb L} \def\bbM{\mathbb M} \def\bbN{\mathbb N} \def\bbO{\mathbb O}
 \def\bbP{\mathbb P} \def\bbQ{\mathbb Q} \def\bbS{\mathbb S} \def\bbT{\mathbb T} \def\bbU{\mathbb U} \def\bbV{\mathbb V}
 \def\bbW{\mathbb W} \def\bbX{\mathbb X} \def\bbY{\mathbb Y} \def\kappa{\varkappa} \def\epsilon{\varepsilon}
 \def\phi{\varphi} \def\le{\leqslant} \def\ge{\geqslant}

\def\UU{\bbU}
\def\Mat{\mathrm{Mat}}
\def\tto{\rightrightarrows}

\def\Gr{\mathrm{Gr}}

\def\graph{\mathrm{graph}}

\def\O{\mathrm{O}}

\def\la{\langle}
\def\ra{\rangle}

\def\B{\mathrm B}
\def\Int{\mathrm {Int}}

\begin{center}
{\bf\Large
 Hua type beta-integrals and projective systems

 \vspace {7pt}
 
 of measures on flag spaces }
 
  \vspace {7pt}

 \sc\large
 Yury A. Neretin%
 \footnote{supported by grant FWF, project P25142}
 
 \end{center}
 
{\small We construct a family of measures on  flag spaces
(or, equivalently, on the spaces of upper-triangular matrices) compatible with respect
to  natural projections. We obtain an $n(n-1)/2$-parametric family of beta-integrals over
space of upper-triangular matrices of size $n$.
}

 \section{Formula}
 
 \COUNTERS
 
 {\bf\punct Beta-integrals.} 
 Let $\K$ be  $\R$,  $\C$, or the quaternions $\H$, let $\kappa=\dim \K$
 (over $\R$). 
 By $dz$ we denote the Lebesgue measure on $\K$. 
 
 Denote by $\cZ_n(\K)$ the space of strictly 
 upper triangular matrices $Z=\{z_{km}\}$ of size $n$, i.e. 
 $ z_{km} =0$ for $k>m$, $z_{kk}=1$, other matrix elements are arbitrary.
 We write $Z^{\{n\}}$ if we wish to emphasize the order of a matrix.
 Denote by $dZ$ or $=dZ^{\{n\}}$ the Lebesgue measure on $\cZ_n(\K)$,
 $$
 dZ=dZ^{\{n\}}:=dz_{12}\dots dz_{1n}\,\,dz_{23}\dots dz_{2n} \,\dots\,dz_{(n-1)n}.
 $$
 The space $\cZ_n(\K)$ is identified in the usual way 
 with the space of flags in $\K^n$ (up to a subset of measure 0). 
 
 For a matrix $Z\in \cZ_n(\K)$ we denote by $[Z]_{pq}$ its left upper corner of size $p\times q$ 
 (it has $p$ rows and $q$ columns). We consider such matrices only for $p<q$. Denote
 $$
 s_{pq}(Z):=\det([Z]_{pq}\,[Z]_{pq}^*).
 $$
 Notice that the matrix $[Z]_{pq}\,[Z]_{pq}^*$ is positive definite and therefore 
 $s_{pq}(Z)$
 are positive.

\begin{theorem} 
Let integer parameters $p$, $q$ range in the domain $1\le p<q\le n$.
For $\lambda_{pq}\in\C$ set
$$
\nu_{pq}:=-\frac12(q-p-1)\kappa+\sum_{k,m:\, p\le k < q,\, q\le m\le n} \lambda_{mk}
.
$$
Then
\begin{equation}
\int_{\cZ_n(\K)} \prod_{1\le p<q\le n} s_{pq}(Z)^{-\lambda_{pq}}\,
dZ^{\{ n\}}=
\pi^{n(n-1)/4} \prod_{1\le p<q\le n} \frac{\Gamma(\nu_{pq}-\kappa/2)}{\Gamma(\nu_{pq})}
\label{eq:main-1}
,\end{equation}
the integral absolutely converges if and only if
$$
\Re \nu_{pq}>\frac 12
.
$$
\end{theorem}

{\sc Remark.}
For $\K=\C$ and $\K=\H$, i.e., $\kappa=2$, $4$, we have cancellations in the right hand side of 
(\ref{eq:main-1}), because
$$\quad\qquad\qquad
\frac{\Gamma(\nu-1)}{\Gamma(\nu)}=\frac 1{\nu-1},\qquad
\frac{\Gamma(\nu-2)}{\Gamma(\nu)}=\frac 1{(\nu-1)(\nu-2)}.\quad\qquad\qquad\boxtimes
$$

Theorem 1.1 is a  corollary of the following statement. 

\begin{theorem}
 Consider the map $\cZ_n(\K)\to \cZ_{n-1}(\K)$ forgetting the last column of a matrix $Z^{\{n\}}$.
 Consider  a measure%
 \footnote{This is a positive measure on $\cZ_n(\K)$ if $\lambda_p\in\R$
 and a complex-valued measure if 
 $\lambda_p\in\C$.}
 \begin{equation}
 \prod\nolimits_{p=1}^{n-1} s_{pn} (z)^{-\lambda_p}\,dZ^{\{n\}}
 \label{eq:density}
 \end{equation}
 on $\cZ_n(\K)$. Assume
 $$
 \lambda_{p}+\lambda_{p+1}+\dots +\lambda_{n-1}>\frac 12(n-p)\kappa
 $$
 for all $p$.
 Then the pushforward of this measure under the forgetting map
 is 
 \begin{multline}
 \label{eq:projection}
\pi^{\frac{(n-1)\kappa}2}\prod_{1\le p\le n-1}
\frac{\Gamma(\lambda_p+\dots +\lambda_n-(n-p)\kappa/2)}
{\Gamma(\lambda_p+\dots +\lambda_n-(n-p+1)\kappa/2)}
\times\\\times
 \prod_{p=1}^{n-2} s_{p(n-1)} (Z^{\{n-1\}})^{-\lambda_p}\,\,dZ^{\{n-1\}}
 .
\end{multline}
\end{theorem}

{\sc Remark.}
On a geometric language the 'forgetting map' is projection of a flag to
a $(n-1)$-dimensional subspace.\hfill $\boxtimes$

\sm

{\bf\punct Comparision with Hua integrals.} In his famous book \cite{Hua}, Hua Loo Keng
obtained a collection of matrix integrals in the following spirit:
$$
I_n(\alpha)=\!\!\!\!\!\!
\int\limits_{T\in \mathrm{Symm}_{n}(\R)}
\!\!\!\!\!\!
\det(1+T^2)^{-\alpha}\,dT
=
\pi^{\frac{n(n+1)}4}\frac{\Gamma(\alpha-n/2)}{\Gamma(\alpha)}
\prod\limits_{p=1}^{n-1}\frac{\Gamma(2\alpha-(n+p)/2)}{\Gamma(2\alpha)}
,$$
where the integration is taken over the space $\mathrm{Symm}_{n}(\R)$
of real symmetric matrices of size $n$.
Recall that $\mathrm{Symm}_{n}(\R)$ is a chart on the real Lagrangian Grassmannian
$\cL_n$ (see, e.g.,
\cite{Ner-gauss}, Sect.3.1, 3.3). Let $\begin{pmatrix}
                                        a&b\\c&d
                                       \end{pmatrix}$ be  a real symplectic matrix.
 The corresponding transformation of $\cL_n$ has the form
 $$\wt T=(a+Tc)^{-1}(b+Td).$$
 The maximal compact subgroup of the real symplectic group
 $\Sp(2n,\R)$
 is isomorphic to $\U(n)$, it can be realized as the group of all real matrices 
 $\begin{pmatrix}
 a&b\\-b&a
  \end{pmatrix}$ such that $a+ib\in\U(n)$. The integrand in the Hua integral
  has nice behavior
  under the action of the group $\U(n)$,
  \begin{equation}
  \det(1+\wt T^2)^{-\alpha}\,d\wt T=\det(1+T^2)^{-\alpha} |\det (a-Tb)|^{-2\alpha +n+1}dT
  \label{eq:hua-inv}
  \end{equation}
  
 Later (see \cite{Gin}, \cite{Ner-beta}) it became known that Hua integrals can be included
 to multi-parametric families of matrix beta-integrals. There is another
  nice property of Hua integrals:
 the measures 
 $$I_n(\sigma+n/2)^{-1}\det(1+T^2)^{-\sigma-n/2} dT$$
 form a projective system under
 the natural 'cutting' map 
 $$\mathrm{Symm}_{n}(\R)\to \mathrm{Symm}_{n-1}(\R).$$
 This phenomenon was a standpoint for infinite-dimensional harmonic analysis
 (see \cite{Pick}, \cite{Ner-hua}, \cite{BO}).
 
 Let us return to our expressions (\ref{eq:density}). To be definite, set $\K=\C$.
 The group $\GL(n,\C)$ acts on the flag space, an element $g\in\GL(n,\C)$
 induces a rational
 transformation  of $\cZ_n(\C)$. For $g\in\GL(n,\C)$
 denote by $\begin{pmatrix}  
               a_p&b_p\\c_p&d_p
              \end{pmatrix}=g$
 its representation as a block matrix of size $p+(n-p)$.
 Denote by $Z=\begin{pmatrix}
               u_p&v_p\\0&w_p
              \end{pmatrix}=Z$
              the corresponding splitting of  $Z\in\cZ_n(\C)$.
 Let $g$ be unitary.
 Then the image of the measure
 (\ref{eq:density}) under the corresponding transformation is
 $$
  \prod_{p=1}^{n-1} \det(u_pa_p+v_pc_p)^{-\lambda_p-4}\cdot
  \prod_{p=1}^{n-1} s_{pn} (Z)^{-\lambda_p}\,dZ
 $$
This is an analog of (\ref{eq:hua-inv}) for flag spaces.

\section{Calculations}

\COUNTERS

{\bf \punct The main lemma.}
Fix $p$ and consider the corner $[z]_{pn}$. Theorem 1.2 is based on the following formula:

\begin{lemma}
\label{l:}
 $$\int_{\K}
 s_{pn}(Z)^{-\lambda} \,dz_{pn}=
 \pi^{\kappa/2} s_{(p-1)(n-1)}^{\lambda-\kappa}(Z)\cdot s_{(p-1)n}^{\kappa/2-\lambda}(Z)
 \cdot s_{p(n-1)}^{\kappa/2-\lambda}(Z).
 $$
\end{lemma}

This lemma is a corollary of two following lemmas.

\begin{lemma}
 Let $a$, $c>0$,  $b\in \K$, and $ac-b\ov b>0$. Let $\lambda\in\C$, $\Re\lambda>\kappa/2$. Then
 $$
 \int_{\K}(a u\ov u+b\ov u+\ov u b+c)^{-\lambda} du=
 \pi^{\kappa/2} a^{\lambda-\kappa}(ac-b\ov b)^{\kappa/2-\lambda}
 \frac {\Gamma(\lambda-\kappa/2)} {\Gamma(\lambda)} .
 $$
\end{lemma}

{\sc Proof.} By an affine change of variable we reduce the integral to the form
$\int_\K (|w|^2+1)^{-\lambda}\,dw$.
\hfill $\square$

\begin{lemma}
\label{l:long}
Denote $u=z_{pn}$. Then  $s_{pn}(Z)$ has the form
 \begin{equation}
 s_{pn}(Z)= a u\ov u+u\ov b+\ov b u+c ,
 \label{eq:trinom}
 \end{equation}
 moreover
 $$
 a=s_{(p-1)(n-1)}(Z), \qquad ac-b\ov b=s_{(p-1)n}(Z)\cdot s_{p(n-1)}(Z)
 .
 $$
\end{lemma}

 Proof of this lemma occupies Subsections 2.2-2.4.
 
 \sm

{\bf\punct Proof of Lemma \ref{l:long} for $\K=\R$ and $\C$.}

{\it Step 0}.
We write $[Z]_{pn}$ as a block matrix of size $\bigl((p-1)+1\bigr)\times \bigl((p-1)+(n-p)+1)\bigr)$,
$$
[Z]_{pn}=\begin{pmatrix}Q&R&t\\ 0&\rho& u  \end{pmatrix}.
$$
Therefore
$$
s_{pn}(Z)=\begin{pmatrix}
               QQ^*+RR^*+tt^*& R\rho^*+t\ov u\\ 
               \rho R^*+u t^*& \rho\rho^*+u\ov u
               \end{pmatrix}
.$$

{\it Step 1.} We write coefficients $a$, $b$, $c$ in (\ref{eq:trinom}):
\begin{multline*}
a=\det(QQ^*+RR^*+tt^*)+\det
  \begin{pmatrix} 
   QQ^*+RR^*+tt^*&t\\t&0
  \end{pmatrix}
  =\\=
  \det  \begin{pmatrix} 
   QQ^*+RR^*+tt^*&t\\t&1
  \end{pmatrix}=
  \det\biggl[  \begin{pmatrix} 
   QQ^*+RR^*&0\\t&1
  \end{pmatrix}
  \begin{pmatrix}
   1&0\\t^*&1
  \end{pmatrix}
  \biggr]
  =\\=
  \det( QQ^*+RR^*) =s_{(p-1)(q-1)}(Z);
\end{multline*}

\begin{multline*}
 \ov b=\det \begin{pmatrix}
    QQ^*+RR^*+tt^*&t\\ \rho R^*&0
   \end{pmatrix}
   =\det\biggl[ \begin{pmatrix}
    QQ^*+RR^*&t\\ \rho R^*&0
   \end{pmatrix}
     \begin{pmatrix}
   1&0\\t^*&1
  \end{pmatrix}\biggr]
 =\\=\det \begin{pmatrix}
    QQ^*+RR^*&t\\ \rho R^*&0
   \end{pmatrix};
\end{multline*}

\begin{equation}
 c=\det\begin{pmatrix}
        QQ^*+RR^*+tt^*&R^*\rho
        \\ R\rho^*&\rho\rho^* \end{pmatrix}.
\end{equation}

{\it Step 2.} Recall the Desnanot-Jacobi formula (see, e.g., \cite{Kra}). 
Consider a block matrix $S$ of size
$m+1+1$,
$$
S=\begin{pmatrix}U&v_1&v_2\\w_1& x_{11}&x_{12}\\w_2& x_{12}&x_{22} \end{pmatrix}.
$$
Then
\begin{multline*}
\det(U)\det(S)=\\=
\det \begin{pmatrix}U&v_1\\w_1& x_{11}\end{pmatrix}\cdot
\det\begin{pmatrix}U&v_2\\w_2& x_{22} \end{pmatrix}-
\det \begin{pmatrix}U&v_1\\w_2& x_{21} \end{pmatrix}\cdot\det 
\begin{pmatrix}U&v_2\\w_1& x_{12}\end{pmatrix}.
\end{multline*}

We apply this identity to the matrix
$$
H:=\begin{pmatrix}
    QQ^*+RR^*+tt^*& R\rho^*&t\\
    \rho R^*&\rho\rho^*&0\\
    t^*&0&1    
   \end{pmatrix}
$$
and get
\begin{equation}
ac-b\ov b= \det(H)\det(QQ^*+RR^*+tt^*).
\label{eq:des}
\end{equation}
Next, we represent $H$ as
\begin{equation}
H=\begin{pmatrix}
    QQ^*+RR^*& R\rho^*&t\\
    \rho R^*&\rho\rho^*&0\\
    0&0&1    
   \end{pmatrix}
   \begin{pmatrix}
    1&0&0\\
    0&1&0\\
    t^*&0&1    
   \end{pmatrix}
   \label{eq:extended}
   .
\end{equation}
Therefore 
\begin{equation}
\det(H)=\begin{pmatrix}
    QQ^*+RR^*& R\rho^*\\
    \rho R^*&\rho\rho^*\\  
   \end{pmatrix}=s_{p(n-1)}(Z)
    \label{eq:extended1}
    .
   \end{equation}
   On the other hand,  $\det(QQ^*+RR^*+tt^*)=s_{(p-1)n}(Z)$,
   and (\ref{eq:des}) implies the desired statement.
   
   \sm
   
   {\bf\punct Several remarks on quaternionic determinants.}
   
   \sm
   
   1) {\it Definition of a quaternionic determinant}. Consider an $m\times m$ matrix 
$A$ over quaternions.
 It determines an automorphism of an $\H$-module
 $\H^m$ and therefore a $\R$-linear operator $A_\R$  in $\R^{4m}$, it can be easily shown that
 $\det(A_\R)\ge 0$.
 We set 
 $$\qquad\qquad\qquad\qquad\qquad  \det(A):=\sqrt[4]{\det(A_\R)}. 
 \qquad\qquad\qquad  \qquad\qquad $$
By definition, $\det(AB)=\det(A)\det(B)$, for strict upper-triangular (and lower triangular)
matrices the determinant is 1, for diagonal matrices with entries $a_{ii}$ it equals
$\prod |a_{ii}|$. These remarks estanblish a coincidence of $\det(A)$ with Dieudonne determinant
\cite{Die} over $\H$.

\sm
   
   2) {\it Formula for determinant of block matrix} 
   (see, e.g., \cite{Gan}, Sect. II.5) remains valid  for quaternionic matrices 
   \begin{equation}
    \det\begin{pmatrix}
         u&v\\w&x
        \end{pmatrix}=\det(u)\cdot \det(x-wu^{-1}v)
        ,
        \label{eq:det-block}
   \end{equation}
here $\begin{pmatrix}
         u&v\\w&x
        \end{pmatrix}$ is a block matrix of size $(m+k)\times(m+k)$, $u$ is assumed to be invertible.
        The formula follows from the identity
        \begin{equation}
        \begin{pmatrix}
         1&0\\-wu^{-1}&1
        \end{pmatrix}
 \begin{pmatrix}
         u&v\\w&x
        \end{pmatrix}
        \begin{pmatrix}
         1&-u^{-1}v\\0&1
        \end{pmatrix}
=\begin{pmatrix}
  u&0\\0&x-wu^{-1}v
 \end{pmatrix}
 \label{eq:proof}.
        \end{equation}
        
  3) {\it Positive definite matrices.} Let a block matrix $\begin{pmatrix}
         u&v\\w&x
        \end{pmatrix}$ be Hermitian positive definite. Then $x-wu^{-1}v$
        also is positive definite (in this case the matrix
        in the left-hand side of (\ref{eq:proof})
   is positive definite).
   
   \sm
 
 4) 
For positive definite block $(m+1)\times (m+1)$ matrices we can write   
\begin{equation}
    \det\begin{pmatrix}
         u&v\\w&x
        \end{pmatrix}=\det(u)\cdot (x-wu^{-1}v).
        \label{eq:idposdef}
\end{equation}
In particular, for $2\times 2$ positive definite matrices we have
\begin{equation}
\begin{pmatrix}u&v\\ \ov v& x\end{pmatrix}=ux- v\ov v.
\label{eq:det22}
\end{equation}
 
 {\bf\punct Proof of Lemma \ref{l:long} for $\K=\H$.} Step 0 is the same.
 
{\it  Step 1.}
 Since $s_{pn}(Z)$ is positive definite,
 we have
\begin{multline*}
s_{pn}(Z)=
\begin{pmatrix}
               QQ^*+RR^*+tt^*& R\rho^*+t\ov u\\ 
               \rho R^*+u t^*& \rho\rho^*+u\ov u
               \end{pmatrix}
               =\\=
\det(QQ^*+RR^*+tt^*)\times \\\times \Bigl[\rho \rho^*+u\ov u- 
(\rho R^*+u t^*)(QQ^*+RR^*+tt^*)^{-1}(R\rho^*+t\ov u) \Bigr]            
.
\end{multline*}
Expanding the expression in square brackets in variables
$u$, $\ov u$, we get
\begin{align}
a&=\det(QQ^*+RR^*+tt^*)\cdot \bigl[1- t^*(QQ^*+RR^*+tt^*)^{-1}t\bigr];
\label{eq:Ha}
\\
b&=\det(QQ^*+RR^*+tt^*)\cdot \bigl[ - t^*(QQ^*+RR^*+tt^*)^{-1}R\rho^*\bigr];
\label{eq:Hb}
\\
c&=\det(QQ^*+RR^*+tt^*)\cdot\bigl[\rho \rho^*-\rho R^*(QQ^*+RR^*+tt^*)^{-1}R\rho^*\bigr].
\label{eq:Hc}
\end{align}
Notice, that  expression (\ref{eq:Ha}) is real and therefore we can write
$ua\ov u=a u\ov u$.

Next, we transform $a$ as
\begin{equation}
a=\det\begin{pmatrix}QQ^*+RR^*+tt^*&t\\t^*&1\end{pmatrix}
\label{eq:Ha1}
\end{equation}
Indeed, this matrix is positive definite because it equals
$XX^*$ with $X=\begin{pmatrix}Q&P&t\\0&0&1 \end{pmatrix}$.
Therefore we can apply  the transformation (\ref{eq:idposdef})
to (\ref{eq:Ha1})
and get the initial expression (\ref{eq:Ha}) for $a$.

This argumentation remains to be valid for $c$ (we do not need a final  expression),
but generally 
$$
b\ne \det \begin{pmatrix}
           QQ^*+RR^*+tt^*&R\rho^*\\ t^*& 0
          \end{pmatrix}
$$
(the expression in the right hand side is real positive, the expression
(\ref{eq:Hb}) for $b$ is not real).

Transforming $a$ as in $\R$-$\C$-cases, we get
$$
a=s_{(p-1)(n-1)}(Z).
$$

{\it Step 2}. The Desnanot-Jacobi identity does not hold for  matrices over $\H$.
However we can adapt previous reasonigs for quaternionic case.
Again, consider the matrix $H$ given by
(\ref{eq:extended}) and transform its determinant
\begin{multline*}
\det
\begin{pmatrix}
    QQ^*+RR^*+tt^*& R\rho^*&t\\
    \rho R^*&\rho\rho^*&0\\
    t^*&0&1    
   \end{pmatrix}
   = 
   \det (QQ^*+RR^*+tt^*)
   \times\\ \times
   \det\left[\begin{pmatrix}
              \rho\rho^*&0\\
    0&1
             \end{pmatrix}+
        \begin{pmatrix} \rho R^*\\ t^* \end{pmatrix}(QQ^*+RR^*+tt^*)^{-1} 
        \begin{pmatrix}
          R\rho^*&t
        \end{pmatrix}\right]
        =\\=
\det (QQ^*+RR^*+tt^*)\times\\ \times
\det\begin{bmatrix}
 1- t^*(QQ^*+RR^*+tt^*)^{-1}t&  - \rho R^*(QQ^*+RR^*+tt^*)^{-1}t
 \\
 - t^*(QQ^*+RR^*+tt^*)^{-1}R\rho^*&
 \rho \rho^*-\rho R^*(QQ^*+RR^*+tt^*)^{-1}R\rho^*
\end{bmatrix}  .     
   \end{multline*}
   Denote by $\Xi$ the $2\times 2$ matrix in square brackets.
The matrix $H$ is positive definite since
$$
H=\begin{pmatrix}
  Q&R&t\\0&\rho&0\\
  0&0&1
  \end{pmatrix}
\begin{pmatrix}
  Q&R&t\\0&\rho&0\\
  0&0&1
  \end{pmatrix}^*,
$$ 
therefore $\Xi$ is positive definite.
Hence we can evaluate $\det[\Xi]$ by (\ref{eq:det22}). Comparing its matrix elements
with (\ref{eq:Ha})--(\ref{eq:Hc}),
we get
$$
\det(\Xi)=\det(QQ^*+RR^*+tt^*)^{-2}(ac-b\ov b)
$$
Therefore
$$
\det (H)= (ac-b\ov b)\cdot \det(QQ^*+RR^*+tt^*)^{-1}=(ac-b\ov b)\cdot s_{(p-1)n}(Z)^{-1}.
$$
The reasoning (\ref{eq:extended})--(\ref{eq:extended1})
remains valid and we get $\det H=s_{p(n-1)}(Z)$. This completes
the calculation.

\sm

 {\bf\punct Proof of Theorem 1.2.}
 Consider a function $\Psi$ on $\cZ_n(\K)$ that does not depend  on variables $z_{1n}$, \dots,
 $z_{(n-1)n}$. We denote such functions as $\Psi(Z^{\{n-1\}})$. Let us transform the integral
 \begin{multline*}
 \int\limits_{\cZ_n(\H)}
 \Psi(Z^{\{n-1\}})\cdot 
  \prod_{p=1}^{n-1} s_{pn} (Z)^{-\lambda_p}\,dZ^{\{n\}}
  =\\=
   \int\limits_{\cZ_{n-1}(\H)\times \K^{n-1}}
 \Psi(Z^{\{n-1\}})\cdot 
  \prod_{p=1}^{n-1} s_{pn} (Z)^{-\lambda_p}\,dZ^{\{n-1\}}\,
  dz_{1n}\,dz_{2n}\dots dz_{(n-1)n}.
  \end{multline*}
  The variable $z_{(n-1)n}$ is presented only in the factor
  $s_{pn} (z)^{-\lambda_{n-1}}(Z)$. Integrating with respect to $z_{(n-1)n}$
  by Lemma \ref{l:}
  we get
 \begin{multline*}
 \pi^{\kappa/2} \frac{\Gamma(\lambda_{n-1}-\kappa/2)}{\Gamma(\lambda_{n-1})}
 \! \!
 \int\limits_{\cZ_{n-2}(\H)\times \K^{n-1}}
  \! \!
 \left[\Psi(Z^{\{n-1\}})\cdot s_{(p-1)(n-1)}
 (Z^{\{n-1\}})^{\lambda_{n-1}-\kappa}
 \right]
   \times\\\times
  \left[\prod_{p}^{n-3} s_{pn} (Z^{\{n\}})^{-\lambda_p}\right]
\cdot
  s_{(n-2)n}(Z^{\{n\}})^{-\lambda_{n-2}-\lambda_{n-1}+\kappa/2} \,dZ^{\{n-1\}}\,
  dz_{1n}\dots dz_{(n-2)n}.
  \end{multline*}
  Now the variable $z_{(n-2)n}$ is presented only in the factor 
  $$s_{(p-1)n}(Z^{\{n\}})^{-\lambda_{n-2}-\lambda_{n-1}+\kappa/2}.$$
  We again apply Lemma \ref{l:}. Etc.
   Finally, we get the integration of
 $\Psi(Z^{\{n-1\}})$ with respect to the measure (\ref{eq:projection}).

\noindent
 {\tt Math.Dept., University of Vienna,
 
 Oskar-Morgenstern-Platz 1, 1090 Wien;
 
\& Institute for Theoretical and Experimental Physics (Moscow)

\& Mech.Math.Dept., Moscow State University,

e-mail: neretin(at) mccme.ru

URL:www.mat.univie.ac.at/$\sim$neretin
} 
 

\begin{thebibliography}{cc}

\bibitem{BO}
Borodin, A.; Olshanski, G.
{\it Harmonic analysis on the infinite-dimensional unitary group and determinantal
point processes.}
Ann. of Math. (2) 161 (2005), no. 3, 1319-1422

\bibitem{Die}
Dieudonn\'e, J.
{\it Les d\'eterminants sur un corps non commutatif.} (French)
Bull. Soc. Math. France 71, (1943). 27-45. 

\bibitem{Gan}
Gantmacher, F. R.
{\it The theory of matrices}.
Reprint of the 1959 translation. AMS Chelsea Publishing, Providence, RI, 1998.

\bibitem{Gin}
Gindikin, S. G.
{\it Analysis in homogeneous domains.} (Russian)
Uspehi Mat. Nauk 19 1964 no. 4 (118), 3-92;
English. transl.: Russ. Math. Surv., 1964, 19:4, 1=89

\bibitem{Hua}
Hua, L. K. {\it Harmonic analysis of functions of several complex variables 
in the classical domains.} Science Press, Peking, 1958 (Chinese); 
Russian translation:
Izdat. Inostr. Lit., Moscow 1959; English translation:
Amer. Math. Soc., Providence, R.I. 1963

\bibitem{Kra}
Krattenthaler, C. {\it Advanced determinant calculus.} The Andrews Festschrift (Maratea, 1998).
S\'em.n  Lothar. Combin. 42 (1999), Art. B42q, 67 pp. (electronic).

\bibitem{Ner-beta}
 Neretin, Yu. A. {\it Matrix analogues of the B-function, and the Plancherel formula 
 for Berezin kernel representations.}  Sb. Math. 191 (2000), no. 5-6, 683-715 

\bibitem{Ner-hua}
Neretin, Yu. A. {\it Hua-type integrals over unitary groups and 
over projective limits of unitary groups.}
Duke Math. J. 114 (2002), no. 2, 239-266.


 
 \bibitem{Ner-gauss}
 Neretin, Yu. A. {\it Lectures on Gaussian integral operators and classical groups.}
 EMS, Z\"urich, 2011.


\bibitem{Pick}
 Pickrell, D.
 {\it Measures on infinite-dimensional Grassmann manifolds.}
 J. Funct. Anal. 70 (1987), no. 2, 323-356.
 
\end{thebibliography}
\end{document}